\documentclass{article}
\usepackage{amsmath,amssymb,amsthm,fullpage}

\providecommand{\abs}[1]{\left\lvert#1\right\rvert}

\DeclareMathOperator{\Aut}{Aut}
\DeclareMathOperator{\Inn}{Inn}
\DeclareMathOperator{\Out}{Out}

\DeclareMathOperator{\Type}{Type}
\DeclareMathOperator{\Core}{Core}

\newtheorem{thm}{Theorem}[section]
\newtheorem{prop}[thm]{Proposition}
\newtheorem{lem}[thm]{Lemma}
\newtheorem{cor}[thm]{Corollary}

\theoremstyle{definition}
\newtheorem{dfn}[thm]{Definition}
\newtheorem*{notn}{Notation}
\newtheorem*{qst}{Question}

\theoremstyle{remark}
\newtheorem*{note}{Remark}

\title{Distance in the Ellipticity Graph}
\author{Yakov Berchenko-Kogan}

\begin{document} 

\maketitle

\begin{abstract}
The ellipticity graph of a free group $F$ was defined by I.~Kapovich and M.~Lustig in order to study the outer automorphism group of $F$, which acts on this graph. The graph was constructed to be analogous to the curve complex of a surface. It is a bipartite graph, whose vertices are conjugacy classes of nontrivial elements of $F$ and equivalence classes of proper free product decompositions of the form $F=A*B$. A conjugacy class is joined by an edge to a free product decomposition $A*B$ whenever the conjugacy class has a representative in either $A$ or $B$. This paper uses Stallings subgroup $X$-digraphs and Whitehead automorphisms to construct algorithms that determine when the distance between two vertices of the ellipticity graph is two, for both types of vertices.
\end{abstract}

\section{Introduction}
A common way to study the outer automorphism group of a free group $F$ is to find geometric spaces on which the group acts. To construct these spaces, it has proven helpful to consider an analogy between the outer automorphism group $\Out(F)$ and another kind of group in geometric group theory, the mapping class group of a surface. For example, the mapping class group acts naturally on a space called Teichm\"uller space. In a key 1986 paper \cite{cv}, Culler and Vogtmann constructed an analogue of Teichm\"uller space which admits an action of the outer automorphism group $\Out(F)$. This space is now known as Outer space. More recently, several papers have turned their attention to another space, called the curve complex, which also admits a natural action of the mapping class group. There appear to be several ways of constructing an analogue of the curve complex that admits an action of $\Out(F)$. Various such analogues are discussed in papers by I.~Kapovich and Lustig \cite{kl}, Behrstock, Bestvina, and Clay \cite{bbc}, Aramayona and Souto \cite{as}, and Sabalka and Savchuk \cite{ss}. One of the analogues of the curve complex defined in \cite{kl} is the ellipticity graph, which we examine more closely in this paper.

\begin{dfn}\label{complex}
Let $F$ be a free group. Given a free product decomposition $A*B$ of $F$, let $[A*B]$ denote the equivalence class of $A*B$ up to conjugating both factors simultaneously and swapping the two factors. The \emph{ellipticity graph} $\mathcal Z(F)$ is a bipartite graph whose vertices are
\begin{equation*}
\{[A*B]\mid A*B\text{ is a proper free product decomposition of }F\}\sqcup\{w\mid w\text{ is a nontrivial conjugacy class of }F\}.
\end{equation*}
The vertices $[A*B]$ and $w$ are adjacent whenever $w$ has a representative in $A$ or in $B$. If this is the case, we say that $w$ is \emph{elliptic} to $A*B$.
\end{dfn}

In Sections \ref{free2} and \ref{cyclic2}, we use subgroup $X$-digraphs and Whitehead automorphisms to construct algorithms that prove the following:

\begin{thm}\label{algsplit}
Let $F$ be a finitely generated free group. There is an algorithm that, given two free product decompositions of $F$, decides whether or not there is a nontrivial element of $F$ that is elliptic to both of them. 
\end{thm}

\begin{thm}\label{algcyc}
Let $F$ be a finitely generated free group. There is an algorithm that, given two conjugacy classes of $F$, decides whether or not they are both elliptic to some proper free product decomposition of $F$.
\end{thm}

Together, these two algorithms can decide whether or not two vertices of the ellipticity graph are a distance of two apart. In addition, one can note that the algorithm of Theorem \ref{algsplit} also decides whether two free product decompositions are adjacent in the dual free splitting graph defined in \cite{kl}, as adjacency in this graph is equivalent to being distance two in the ellipticity graph.

During the review process, I was informed that a statement equivalent to Theorem \ref{algcyc} appears as the main result in a 1996 paper by Stallings \cite{swgoh}, with a proof that uses similar techniques. In addition, Theorem \ref{algsplit} follows from Propositions 9.7 and 9.8 in a paper by I.~Kapovich and Myasnikov \cite{km}. As such, the results of this paper are not new, but I hope that this exposition will be found useful.

\section{Definitions and Supporting Results}
\subsection{The ellipticity graph}\label{free}
The vertices of the ellipticity graph are conjugacy classes of elements and equivalence classes of free product decompositions. We begin by defining the notation we will use for the free group and for cyclic words.

\begin{dfn}
If $X$ is a set, let $\Sigma=X\cup X^{-1}$. We denote the \emph{free group on $X$} by $F(X)$ and interpret it as the set of all freely reduced words in $\Sigma$. The \emph{length} of an element $g$ of $F(X)$ is denoted $\abs{g}_X$ or simply $\abs g$.
\end{dfn}

\begin{dfn}
We identify each conjugacy class of $F(X)$ with a \emph{cyclic word}, defined as a cyclically ordered set of letters of $\Sigma$ such that no two consecutive letters are inverses of each other. The \emph{length} of a cyclic word $w$ is denoted $\abs{w}_X$ or simply $\abs w$.
\end{dfn}

\begin{dfn}
We say that an element $g$ of the free group $F(X)$ is \emph{cyclically reduced} if the word's first letter is not the inverse of the word's last letter. Equivalently, if $w$ is the cyclic word corresponding to $g$, then $g$ is cyclically reduced if $\abs g=\abs w$.
\end{dfn}

Next we define free product decompositions, as well as their equivalence classes.
\begin{dfn}[free product decompositions]
Let $G$ be a group. A \emph{free product decomposition} of $G$ is a decomposition of $G$ as a free product of subgroups $A*B$. We say that such a splitting is proper if both factors are proper subgroups.
\end{dfn}

\begin{dfn}[Equivalence classes of free product decompositions]\label{equiv}
If $F$ is a free group, we define equivalence classes of free product decompositions of $F$ up to conjugation and swapping the factors. Namely, we let the equivalence class $[A*B]$ contain all $C*D$ such that either $A=xCx^{-1}$ and $B=xDx^{-1}$ or $A=xDx^{-1}$ and $B=xCx^{-1}$ for some $x\in F$.
\end{dfn}

\begin{dfn}[Ellipticity]
Let $A*B$ be a free product decomposition. An element $g\in G$ is \emph{elliptic} to $A*B$ if it is conjugate either to an element of $A$ or to an element of $B$. This definition naturally extends to the notion of ellipticity of conjugacy classes of elements to equivalence classes of free product decompositions in Definition \ref{complex}.
\end{dfn}

\begin{note}
The group $\Aut(F)$ acts on the ellipticity graph $\mathcal Z(F)$ preserving edges. Since the group of inner automorphisms $\Inn(F)$ fixes $\mathcal Z(F)$, the group $\Out(F)=\Aut(F)/\Inn(F)$ acts on $\mathcal Z(F)$.
\end{note}

The ellipticity graph is not locally finite, so computing distances is not easy. Nevertheless, some questions about distance in the graph can be answered.

\subsection{$X$-digraphs}\label{digraphs}

Stallings digraphs are an elegant way of representing a subgroup of a free group as a graph with labeled edges, first discussed in a paper by John Stallings \cite{s}. One can use these graphs to produce simple algorithms for various problems, such as deciding whether a given word is in the subgroup. We present the relevant results here. For a more complete treatment of Stallings subgroup digraphs, see the paper by I.~Kapovich and Myasnikov \cite{km}.

\begin{dfn}
If $X$ is a finite alphabet, then an \emph{$X$-labeled digraph} or \emph{$X$-digraph} $\Gamma$ is a 5-tuple $(V,E,o,t,\mu)$, where $V$ and $E$ are sets and $o$, $t$, and $\mu$ are functions with $o,t\colon E\rightarrow V$ and $\mu\colon E\rightarrow X$. The set $V$, also denoted $V\Gamma$, contains the \emph{vertices} of $\Gamma$. The set $E$, also denoted $E\Gamma$, contains the \emph{edges} of $\Gamma$. For any edge $e$, $o(e)$, $t(e)$, and $\mu(e)$ are the \emph{origin}, \emph{terminus}, and \emph{label} of $e$, respectively.

A \emph{morphism} between two $X$-digraphs is a map that sends vertices to vertices, sends edges to edges, sends the origin and terminus of an edge to the origin and terminus of the image of that edge, and preserves the labels of edges.
\end{dfn}

\begin{dfn}
An \emph{$X$-digraph with base vertex} is a pair $(\Gamma,v)$ where $\Gamma$ is an $X$-digraph and $v$ is a vertex of $\Gamma$. A \emph{morphism of $X$-digraphs with base vertex} is a morphism of $X$-digraphs that sends the base vertex of one $X$-digraph to the base vertex of the other.
\end{dfn}

There is a bijection between certain kinds of $X$-digraphs with base vertex and subgroups of $F(X)$. In order to define this bijection, we first define paths.

\begin{dfn}
Let $\Sigma=X\cup X^{-1}$. To each $X$-digraph $\Gamma$ we will associate a $\Sigma$-digraph $\hat\Gamma$ with vertex set $\hat V=V$ and edge set $\hat E=E\cup E^{-1}$. For each edge $e$ of $\Gamma$ from $v$ to $w$ with label $x$, the edge $e^{-1}$ of $\hat\Gamma$ is defined to go from $w$ to $v$ with label $x^{-1}$.

Every edge $\hat e$ of $\hat\Gamma$ has a corresponding edge of $\Gamma$, which we will call the \emph{positive edge corresponding to $\hat e$}.
\end{dfn}

\begin{note}
It is clear that we can recover $\Gamma$ from $\hat\Gamma$ by taking the positive edges of $\hat\Gamma$ and restricting $o$, $t$, and $\mu$ appropriately. Therefore, we will occasionally abuse notation and refer to $\hat\Gamma$ as $\Gamma$.
\end{note}

\begin{dfn}
A \emph{nontrivial path} $p$ in $\Gamma$ is a finite sequence of edges of $\hat\Gamma$ such that the terminus of each edge is the origin of the next edge. We define $o(p)$ and $t(p)$ naturally. The label $\mu(p)$ is the word over the alphabet $\Sigma$ constructed by writing the labels of the edges of $p$ in order. Note that $\mu(p)$ is not necessarily freely reduced. We let $\overline{\mu(p)}$ denote the corresponding freely reduced word.

A \emph{trivial path} $p$ is a vertex $v$ of $\Gamma$. We have $o(p)=t(p)=v$, and the label of $p$ is the empty word.
\end{dfn}

\begin{notn}
If $p$ is a path of $\Gamma$ with $p=e_1e_2\cdots e_k$, then let $\tilde p$ denote the subgraph of $\Gamma$ containing all positive edges corresponding to edges of $p$, as well as the origins and termini of all those edges. If $p$ is a trivial path of $\Gamma$ from a vertex $v$ to itself, let $\tilde p$ denote the subgraph of $\Gamma$ containing only the vertex $v$ and no edges.
\end{notn}

In general, we will only concern ourselves with paths that do not backtrack.

\begin{dfn}
A path $p$ in an $X$-digraph $\Gamma$ is \emph{reduced} if it does not contain an edge $e\in E\hat\Gamma$ such that the next edge in the path is $e^{-1}$. A \emph{path reduction} of a nonreduced path $p$ is a construction of a path $p'$ by removing of two consecutive edges of $p$ that are inverses of each other. To \emph{reduce} a path $p$ means to construct a path $\bar{p}$ by performing all possible path reductions.
\end{dfn}

\begin{note}
One can check that $o(\bar p)=o(p)$, $t(\bar p)=t(p)$, and $\overline{\mu(\bar p)}=\overline{\mu(p)}$.
\end{note}

We use paths to define the language of an $X$-digraph with base vertex.

\begin{dfn}
If $\Gamma$ is an $X$-digraph and $v$ is a vertex of $\Gamma$, then the \emph{language of $\Gamma$ with respect to $v$}, denoted $L(\Gamma,v)$, is the set of all labels of reduced paths from $v$ to $v$. That is,
\begin{equation*}
L(\Gamma,v)=\{\mu(p)\mid p\text{ is a reduced path in }\Gamma\text{ and }o(p)=t(p)=v\}.
\end{equation*}
\end{dfn}

Usually, we would like our $X$-digraphs to satisfy two additional properties: We would like them to be \emph{folded} and \emph{core}.

\begin{dfn}
An $X$-digraph $\Gamma$ is \emph{folded} if for every vertex $v$ and label $x$, there is at most one edge with origin $v$ and label $x$, and there is at most one edge with terminus $v$ and label $x$.
\end{dfn}

\begin{note}
Any finite $X$-digraph can be easily transformed into a folded $X$-digraph in a way that preserves $\overline{L(\Gamma,v)}$, the set of free reductions of the words in the language of the digraph. This transformation is called \emph{Stallings folding} and is discussed in \cite{km}.
\end{note}

\begin{dfn}
Let $\Gamma$ be an $X$-digraph and $v$ be a vertex of $\Gamma$. The \emph{core of $\Gamma$ at $v$}, denoted by $\Core(\Gamma,v)$, is the union of all of the reduced paths from $v$ back to $v$. That is, $\Core(\Gamma,v)$ is a subgraph of $\Gamma$ with base vertex $v$ defined by
\begin{equation*}
\Core(\Gamma,v):=\bigcup\{\tilde p\mid p\text{ is a reduced path in $\Gamma$ from $v$ to $v$}\}.
\end{equation*}

If $\Core(\Gamma,v)=(\Gamma,v)$, we will say that $(\Gamma,v)$ is \emph{core} or $\Gamma$ is a \emph{core graph with respect to $v$}.
\end{dfn}

\begin{note}
Since $\Core(\Gamma,v)$ and $(\Gamma,v)$ have the same reduced paths, we have $L(\Core(\Gamma,v))=L(\Gamma,v)$.
\end{note}

The following theorems show that there is a bijection between finitely generated subgroups of $F(X)$ and finite, folded, and core $X$-digraphs.

\begin{thm}[{\cite[Lemma 3.2]{km}}]
If $(\Gamma,v)$ is a folded $X$-digraph with base vertex, then $L(\Gamma,v)$ is a subgroup of $F(X)$.
\end{thm}

\begin{thm}[{\cite[Propositions 3.8, 5.1, and 5.2]{km}}]
For every subgroup $H$ of the free group $F(X)$ there is an $X$-digraph with base vertex, denoted $(\Gamma(H),1_H)$, such that $\Gamma(H)$ is folded and core with respect to $1_H$ and such that $L(\Gamma(H),1_H)=H$. Moreover, $(\Gamma(H),1_H)$ is unique up to isomorphism of $X$-digraphs with base vertex.

If $H$ is finitely generated, then $\Gamma(H)$ is a finite graph that can be constructed in finite time.
\end{thm}

\begin{cor}\label{glg}
If $(\Gamma,v)$ is folded and core, then
\begin{equation*}
\Bigl(\Gamma\bigl(L(\Gamma,v)\bigr),1_{L(\Gamma,v)}\Bigr)\cong(\Gamma,v).
\end{equation*}
\begin{proof}
Both $\bigl(\Gamma\bigl(L(\Gamma,v)\bigr),1_{L(\Gamma,v)}\bigr)$ and $(\Gamma,v)$ are folded, core, and have language $L(\Gamma,v)$.
\end{proof}
\end{cor}

We now present the applications of $X$-digraphs that are relevant to this paper. First, we describe an algorithm for determining when two subgroups are conjugate to each other.

\begin{dfn}\label{type}
If $\Gamma$ is core with respect to some vertex $v$, we will define the \emph{type} of $\Gamma$, denoted $\Type(\Gamma)$, as follows.

If $v$ does not have degree one in $\Gamma$, then $\Type(\Gamma)$ is defined to be $\Gamma$.

Otherwise, let $p$ be the unique nontrivial reduced path originating at $v$ such that vertices of $\tilde p$ other than $o(p)$ and $t(p)$ have degree two in $\Gamma$, and $t(p)$ has degree greater than two. Then $\Type(\Gamma)$ is defined to be the subgraph of $\Gamma$ constructed by removing from $\Gamma$ all edges and all vertices of $\tilde p$, except for $t(p)$.

Equivalently, if $\Gamma$ is core with respect to some vertex $v$, then
\begin{equation*}
\Type(\Gamma):=\bigcap_{u\in V\Gamma}\Core(\Gamma,u).
\end{equation*}
\end{dfn}

\begin{thm}[{\cite[Proposition 7.7]{km}}]\label{conj}
If $H$ and $K$ are subgroups of $F(X)$, then $H$ is conjugate to $K$ if and only if $\Type(\Gamma(H))$ and $\Type(\Gamma(K))$ are isomorphic as $X$-digraphs.
\end{thm}

Next, we present the construction of the digraph corresponding to the intersection of two subgroups.

\begin{dfn}
If $\Gamma$ and $\Delta$ are $X$-digraphs, then the \emph{product} graph, denoted $\Gamma\times\Delta$, is defined as follows.
\begin{align*}
V(\Gamma\times\Delta)&:=V\Gamma\times V\Delta,\\
E(\Gamma\times\Delta)&:=\{(e,f)\in E\Gamma\times E\Delta\mid\mu(e)=\mu(f)\},\\
o(e,f)&:=(o(e),o(f)),\\
t(e,f)&:=(t(e),t(f)),\\
\mu(e,f)&:=\mu(e)=\mu(f).
\end{align*}
\end{dfn}

\begin{note}
The product graph $\Gamma\times\Delta$ is not necessarily connected, even if both $\Gamma$ and $\Delta$ are connected.
\end{note}

\begin{thm}[{\cite[Lemma 9.3]{km}}]\label{prodl}
If $(\Gamma,u)$ and $(\Delta,v)$ are two folded $X$-digraphs with base vertex, then $\Gamma\times\Delta$ is folded and
\begin{equation*}
L\bigl(\Gamma\times\Delta,(u,v)\bigr)=L(\Gamma,u)\cap L(\Delta,v).
\end{equation*}
\end{thm}

\begin{thm}[{\cite[Proposition 9.4]{km}}]\label{prod}
If $H$ and $K$ are two subgroups of $F(X)$, then
\begin{equation*}
\bigl(\Gamma(H\cap K),1_{H\cap K}\bigr)\cong\Core\bigl(\Gamma(H)\times\Gamma(K),(1_H,1_K)\bigr).
\end{equation*}
\end{thm}

\subsection{Whitehead automorphisms}\label{nielwhit}
The Whitehead automorphisms are a set of generators for $\Aut(F(X))$ that can be used to determine when tuples of cyclic words are in the same orbit.

\begin{dfn}
A \emph{Whitehead automorphism} $\tau$ of $F(X)$ is an automorphism of $F(X)$ satisfying one of the following two properties.
\begin{itemize}
\item The automorphism $\tau$ permutes the elements of the set $X\cup X^{-1}$. In this we case will call $\tau$ a \emph{relabeling automorphism}.

\item The set $X\cup X^{-1}$ contains a letter $a$, called the \emph{multiplier} of $\tau$, such that for all $x\in X\cup X^{-1}$, we have $\tau(x)\in\{x,xa,a^{-1}x,a^{-1}xa\}$.
\end{itemize}
The set of all Whitehead automorphisms is denoted $\Omega$.
\end{dfn}

The following theorems describe Whitehead's algorithm.

\begin{thm}[{\cite[Proposition~4.20]{ls}}]\label{white}
Suppose $w_1,\dotsc,w_t$ and $w_1',\dotsc,w_t'$ are cyclic words in $F$ such that for some $\alpha\in\Aut(F)$ we have $\alpha(w_h)=w_h'$ for ${1\le h\le t}$. Suppose that the sum $\sum{\abs{w_h'}}$ is minimal among all sums of the form $\sum{\abs{\alpha'(w_h)}}$ for $\alpha'\in\Aut(F)$. Then there exist Whitehead automorphisms $\tau_1,\dotsc,\tau_n$ such that $\alpha=\tau_n\dotsm\tau_1$ and $\sum{\abs{(\tau_i\dotsm\tau_1)(w_h)}}\le\sum{\abs{w_h}}$ for ${0<i<n}$ with strict inequality unless $\sum{\abs{w_h}}=\sum{\abs{w_h'}}$.
\end{thm}

In other words, if a tuple of cyclic words has minimal length in its orbit, then one can arrive at this tuple from any other tuple in the orbit by a sequence of Whitehead automorphisms that first strictly decrease the length of the tuple until the length is minimal and then keep the length of the tuple the same.

\begin{cor}\label{whitcor}
If the cyclic words $w_1,\dotsc,w_t$ are such that $\sum{\abs{w_h}}$ is not minimal in the orbit of $(w_1,\dotsc,w_t)$ under the action of $\Aut(F)$, then there exists a Whitehead automorphism $\tau$ such that $\sum{\abs{\tau(w_h)}}<\sum{\abs{w_h}}$.
\end{cor}

This theorem is the key ingredient in Whitehead's algorithm:

\begin{thm}[{\cite[Proposition~4.21]{ls}}]\label{orbit}
If $w_1,\dotsc,w_t$ and $w_1',\dotsc,w_t'$ are cyclic words then it is decidable whether or not there exists an automorphism $\alpha$ of $F$ such that $\alpha(w_h)=w_h'$ for $1\le h\le t$.
\end{thm}

\section{Free Product Decompositions with a Common Elliptic Element}\label{free2}
We aim to determine whether or not the distance in the ellipticity graph between two given classes of free product decompositions is two. In order to do this, we first answer the following question.

\begin{qst}
Given a finitely generated free group $F$ and two finitely generated subgroups $H$ and $K$ of $F$, how can we determine whether or not there exists a nontrivial element $g$ of $F$ such that $g$ is conjugate to an element of $H$ and $g$ is conjugate to an element of $K$?
\end{qst}

\begin{notn}
Given two subgroups $H$ and $K$, let $C_{H,K}$, or simply $C$, denote the set
\begin{equation*}
C:=\bigcup_{x,y\in F}xHx^{-1}\cap yKy^{-1}
\end{equation*}
\end{notn}

The set $C$ is the set of all elements of $F$ that are conjugate to an element of $H$ and conjugate to an element of $K$. In particular, if $g\in C$, then every element in the conjugacy class of $g$ is also in $C$.

We aim, given $H$ and $K$, to decide whether or not $C=\{1\}$. We first establish the following lemma.

\begin{lem}
If $\Gamma$ is a folded $X$-digraph with base vertex $v$ and there is a nontrivial cyclically reduced element $g$ of $L(\Gamma,v)$, then $v$ has degree at least two.

\begin{proof}
Let $x$ be the first letter of $g$, and let $y$ be the last letter of $g$. There is a reduced path in $\Gamma$ from $v$ to $v$ such that the label of the first edge $e$ of the path is $x$ and the label of the last edge $f$ of the path is $y$. Note that $e$ and $f^{-1}$ have origin $v$ and labels $x$ and $y^{-1}$, respectively. Since $g$ is cyclically reduced, $x\neq y^{-1}$. Thus $e$ and $f$ cannot be the same edge, so there are at least two edges originating from $v$.
\end{proof}
\end{lem}

\begin{cor}\label{equaltotype}
If $\Gamma$ is folded and a core graph with respect to $v$, and $L(\Gamma,v)$ contains a nontrivial cyclically reduced element, then $\Type(\Gamma)=\Gamma$.

\begin{proof}
This is a straightforward application of Definition \ref{type}.
\end{proof}
\end{cor}

\begin{lem}
There exists a nontrivial element of $C$ if and only if there exists a vertex $u$ of $\Type(\Gamma(H))$ and a vertex $v$ of $\Type(\Gamma(K))$ such that \begin{equation*}
L\bigl(\Type(\Gamma(H))\times \Type(\Gamma(K)),(u,v)\bigr)\neq 1.
\end{equation*}

\begin{proof}
If $g$ is a nontrivial element of $C$, then let $g'$ be a cyclically reduced element of the conjugacy class of $g$. We know that $g'\in C$. Therefore there exist $x$ and $y$ such that $g'\in xHx^{-1}$ and $g'\in yKy^{-1}$. By Theorem \ref{conj}, $\Type(\Gamma(H))$ is isomorphic to $\Type(\Gamma(xHx^{-1}))$, which, in turn, is equal to $\Gamma(xHx^{-1})$ by Corollary \ref{equaltotype}. Similarly, $\Type(\Gamma(K))$ is isomorphic to $\Gamma(yKy^{-1})$. Therefore, there exists a vertex $u$ of $\Type(\Gamma(H))$ and a vertex $v$ of $\Type(\Gamma(K))$ such that
\begin{align*}
(\Type(\Gamma(H)),u)&\cong(\Gamma(xHx^{-1}),1_{xHx^{-1}})\text{ and}\\
(\Type(\Gamma(K)),v)&\cong(\Gamma(yKy^{-1}),1_{yKy^{-1}}).
\end{align*}
Using Theorem \ref{prod}, we have that
\begin{equation*}
\begin{split}
g'&\in xHx^{-1}\cap yKy^{-1}\\
&=L\bigl(\Gamma(xHx^{-1}\cap yKy^{-1}),1_{xHx^{-1}\cap yKy^{-1}}\bigr)\\
&=L\bigl(\Core\bigl(\Gamma(xHx^{-1})\times\Gamma(yKy^{-1}),(1_{xHx^{-1}},1_{yKy^{-1}})\bigr)\bigr)\\
&=L\bigl(\Gamma(xHx^{-1})\times\Gamma(yKy^{-1}),(1_{xHx^{-1}},1_{yKy^{-1}})\bigr)\\
&=L\bigl(\Type(\Gamma(H))\times \Type(\Gamma(K)),(u,v)\bigr)
\end{split}
\end{equation*}
Since $g'$ is nontrivial, $L\bigl(\Type(\Gamma(H))\times \Type(\Gamma(K)),(u,v)\bigr)\neq1$.

Conversely, assume there does exist a choice of $u$ and $v$ such that $L\bigl(\Type(\Gamma(H))\times \Type(\Gamma(K)),(u,v)\bigr)$ contains a nontrivial element $g$. Then by Theorem \ref{prodl}, we have that 
\begin{equation*}
g\in L\bigl(\Type(\Gamma(H)),u\bigr)\cap L\bigl(\Type(\Gamma(K)),v\bigr).
\end{equation*}
By Corollary \ref{glg},
\begin{equation*}
\Type\Bigl(\Gamma\bigl(L\bigl(\Type(\Gamma(H)),u\bigr)\bigr)\Bigr)\cong \Type\bigl(\Type(\Gamma(H))\bigr)=\Type(\Gamma(H)).
\end{equation*}
Therefore by Theorem \ref{conj}, the subgroup $L\bigl(\Type(\Gamma(H)),u\bigr)$ is conjugate to the subgroup $H$. Similarly, $L\bigl(\Type(\Gamma(K)),v\bigr)$ is conjugate to $K$. Thus, $g\in xHx^{-1}\cap yKy^{-1}$ for some $x,y\in F$, so $g$ is a nontrivial element of $C$.
\end{proof}
\end{lem}

\begin{lem}
We have that $L\bigl(\Type(\Gamma(H))\times \Type(\Gamma(K)),(u,v)\bigr)=1$ for all vertices $(u,v)$ if and only if $\Type(\Gamma(H))\times \Type(\Gamma(K))$ is acyclic.

\begin{proof}
If $\Type(\Gamma(H))\times \Type(\Gamma(K))$ is acyclic, then for any vertex $(u,v)$ of $\Type(\Gamma(H))\times \Type(\Gamma(K))$, the only reduced path from $(u,v)$ to $(u,v)$ is the trivial path, so $L\bigl(\Type(\Gamma(H))\times \Type(\Gamma(K)),(u,v)\bigr)=1$. Conversely, assume that $\Type(\Gamma(H))\times \Type(\Gamma(K))$ has a cycle. Let $(u,v)$ be a vertex on this cycle, and let $p$ be a reduced path from $(u,v)$ to $(u,v)$ going once around the cycle. Then $\mu(p)\neq1$ and $\mu(p)\in L\bigl(\Type(\Gamma(H))\times \Type(\Gamma(K)),(u,v)\bigr)$, so $L\bigl(\Type(\Gamma(H))\times \Type(\Gamma(K)),(u,v)\bigr)\neq1$.
\end{proof}
\end{lem}

The above lemmas establish the following two facts.

\begin{prop}
Given a finitely generated free group $F$ and two finitely generated subgroups $H$ and $K$ of $F$, there exists a nontrivial element of $F$ conjugate to both an element of $H$ and an element of $K$ if and only if the graph $\Type(\Gamma(H))\times \Type(\Gamma(K))$ has a cycle.
\end{prop}

Since these graphs are finite, this problem is decidable. We can now prove Theorem \ref{algsplit}.

\begin{thm}
If $F$ is a finitely generated free group, then it is decidable whether or not two classes of proper free product decompositions are distance two in the ellipticity graph.

\begin{proof}
Let $[A*B]$ and $[C*D]$ be the two classes of proper free product decompositions. The two classes are distance two if and only if there exists a nontrivial conjugacy class elliptic to both splittings. Such a conjugacy class exists if and only if there exists a nontrivial element $g$ of $F$ such that $g$ is conjugate to an element of either $A$ or $B$ and conjugate to an element of either $C$ or $D$. Since free factors of finitely generated free groups are finitely generated, we can decide whether or not there exists a nontrivial element $g$ conjugate to both $A$ and $C$. Likewise, we can test if such an element exists for the pairs of groups $A$ and $D$, $B$ and $C$, and $B$ and $D$. Therefore, we can decide whether or not $[A*B]$ and $[C*D]$ are distance two in $\mathcal Z(F)$.
\end{proof}
\end{thm}

\begin{note}
If $A*B$ and $C*D$ are two free product decompositions, we can choose a free basis $X_A$ for $A$ and a free basis $X_B$ for $B$. We can then define a free basis for $F$ by $X:=X_A\cup X_B$. Then $\Gamma(A)$ and $\Gamma(B)$ each contain exactly one vertex, so the product graphs are easy to construct. For example, $\Gamma(A)\times\Gamma(C)$ is isomorphic to the subgraph of $\Gamma(C)$ with all edges with labels in $X_B$ removed.
\end{note}

\section{Elements Elliptic to a Common Free Product Decomposition}\label{cyclic2}
We now turn to the distance two problem for conjugacy classes of elements.

\begin{qst}
If $F$ is a free group with free basis $X$ and $v$ and $w$ are two nontrivial conjugacy classes of $F$, how can we determine whether or not there exists a proper free product decomposition $A*B$ of $F$ such that both $v$ and $w$ are elliptic to $A*B$?
\end{qst}

We identify the conjugacy classes $v$ and $w$ with their corresponding cyclic words in the basis $X$, and we define some shorthand for talking about pairs of such cyclic words. 

\begin{notn}\label{letters}
If $w$ is a cyclic word, let $\mathcal L(w)$ denote the subset of $X$ containing all letters $x$ such that either $x$ or $x^{-1}$ appears in $w$. Essentially, $\mathcal L(w)$ denotes the letters in $X$ used by $w$. Let $\Lambda(w)$ denote $\mathcal L(w)\cup (\mathcal L(w))^{-1}$.
\end{notn}

\begin{dfn}\label{good}
We call a pair $(v,w)$ \emph{frugal} if $v$ and $w$ do not use all letters of $X$, that is, $\mathcal L(v)\cup\mathcal L(w)\neq X$. We call $(v,w)$ \emph{disjoint} if $v$ and $w$ do not share any letters, that is, $\mathcal L(v)\cap\mathcal L(w)=\emptyset$. We call $(v,w)$ \emph{good} if it is frugal or disjoint. Let the \emph{length} of a pair $(v,w)$ be defined to be $\abs v+\abs w$. We say that a pair $(v,w)$ has \emph{minimal length} if its length is minimal in its orbit under $\Aut(F)$. In other words, $(v,w)$ has minimal length if $\abs v+\abs w\le\abs{\phi(v)}+\abs{\phi(w)}$ for all $\phi\in\Aut(F)$.
\end{dfn}

We can reduce the problem of finding a free product decomposition adjacent to two cyclic words to the problem of finding an automorphism of $F$ satisfying certain properties.

\begin{lem}\label{auto}
Let $v$ and $w$ be nontrivial cyclic words. Then there exists a proper free product decomposition $A*B$ of $F$ such that both $v$ and $w$ are elliptic to $A*B$ if and only if there exists an automorphism $\phi$ of $F$ such that $(\phi(v),\phi(w))$ is good.

\begin{proof}
If there exists a proper free product decomposition $A*B$ of $F$ such that both $v$ and $w$ are elliptic to $A*B$, let $\phi$ be an automorphism that sends $A$ to $\langle X_1\rangle$ and $B$ to $\langle X_2\rangle$ for some subsets $X_1$ and $X_2$ of $X$. Since $A*B$ is a free product decomposition, $X_1\sqcup X_2=X$. Since the free product decomposition is proper, we know that $X_1$ and $X_2$ are both proper subsets of $X$. It is clear that $\phi(v)$ and $\phi(w)$ are elliptic to $\phi(A)*\phi(B)=\langle X_1\rangle*\langle X_2\rangle$.

Since $\phi(v)$ is elliptic to $\langle X_1\rangle*\langle X_2\rangle$, it has a representative in either $\langle X_1\rangle$ or $\langle X_2\rangle$. Assume without loss of generality that $\phi(v)$ has a representative $g$ in $\langle X_1\rangle$. Since all the letters of $\phi(v)$ are in the word $g$, we know $\mathcal L(\phi(v))\subseteq X_1$. Similarly, either $\mathcal L(\phi(w))\subseteq X_1$ or $\mathcal L(\phi(w))\subseteq X_2$. If $\mathcal L(\phi(w))\subseteq X_1$, then $(\phi(v),\phi(w))$ is frugal, and if $\mathcal L(\phi(w))\subseteq X_2$, then $(\phi(v),\phi(w))$ is disjoint.

Conversely, if there exists an automorphism $\phi$ of $F$ such that $\mathcal L(\phi(v))\cup\mathcal L(\phi(w))\neq X$, then let $X_1=\mathcal L(\phi(v))\cup\mathcal L(\phi(w))$ and let $X_2=X\setminus X_1$. If, on the other hand, there exists an automorphism $\phi$ of $F$ such that $\mathcal L(\phi(v))\cap\mathcal L(\phi(w))=\emptyset$, then let $X_1=\mathcal L(\phi(v))$ and let $X_2=X\setminus X_1$.

In both cases $\langle X_1\rangle*\langle X_2\rangle$ is a proper free product decomposition of $F$, so let $A=\phi^{-1}(\langle X_1\rangle)$ and $B=\phi^{-1}(\langle X_2\rangle)$. Both $\phi(v)$ and $\phi(w)$ are elliptic to $\langle X_1\rangle*\langle X_2\rangle$, so both $v$ and $w$ are elliptic to $A*B$.
\end{proof}
\end{lem}

In order to determine whether or not there is a good pair in the orbit of $(v,w)$, we prove some results about Whitehead automorphisms.

\begin{lem}\label{badmult}
If $w$ is a cyclic word and $\tau$ is a Whitehead automorphism with multiplier $a\notin\Lambda(w)$, then either $\tau(w)=w$ or $\abs{\tau(w)}>\abs w$.

\begin{proof}
Let $w=w_1w_2\cdots w_l$ where $w_i$ are letters. Then $\tau(w)$ is the cyclic word $w_1a^{\epsilon_1}w_2a^{\epsilon_2}\cdots w_la^{\epsilon_l}$ for some exponents $\epsilon_i$. Since $a\neq w_i\neq a^{-1}$ for all $i$, this cyclic word is reduced. If all of the exponents $\epsilon_i$ are zero, then $\tau(w)=w$. If one of the exponents $\epsilon_i$ is nonzero, then the cyclic word $w_1a^{\epsilon_1}w_2a^{\epsilon_2}\cdots w_la^{\epsilon_l}$ has at least one more letter than the cyclic word $w_1w_2\cdots w_l$, so $\abs{\tau(w)}>\abs w$.
\end{proof}
\end{lem}

\begin{cor}\label{badmultcor}
If $\tau$ is a Whitehead automorphism with multiplier $a$ and $\abs{\tau(w)}\le\abs w$, then $\mathcal L(\tau(w))\subseteq\mathcal L(w)$.
\begin{proof}
If $a\in\Lambda(w)$, then $\mathcal L(\tau(w))\subseteq\mathcal L(w)$. If $a\notin\Lambda(w)$, then by Lemma \ref{badmult}, since it cannot be that $\abs{\tau(w)}>\abs w$, we have that $\tau(w)=w$, so $\mathcal L(\tau(w))=\mathcal L(w)$.
\end{proof}
\end{cor}

\begin{lem}
If there is a good pair $(v,w)$ of cyclic words that does not have minimal length, then there exists a Whitehead automorphism $\tau$ such that $\bigl(\tau(v),\tau(w)\bigr)$ is also a good pair and such that the length of $\bigl(\tau(v),\tau(w)\bigr)$ is smaller than the length of $(v,w)$.

\begin{proof}
Since $(v,w)$ does not have minimal length, by Corollary \ref{whitcor}, there exists a Whitehead automorphism $\tau_1$ such that $\abs{\tau_1(v)}+\abs{\tau_1(w)}<\abs v+\abs w$. Clearly $\tau_1$ is not a relabeling automorphism.

If $(v,w)$ is frugal, then let $\tau=\tau_1$. We will show that $(\tau(v),\tau(w))$ is frugal. Let the letter $a$ be the multiplier of the Whitehead automorphism $\tau$. If $a$ is not in $\Lambda(v)\cup\Lambda(w)$, then by Lemma \ref{badmult}, $\abs{\tau(v)}\ge\abs v$ and $\abs{\tau(w)}\ge\abs w$, which is a contradiction. Thus $a\in\Lambda(v)\cup\Lambda(w)$, so $\mathcal L(\tau(v))\subseteq\mathcal L(v)\cup\mathcal L(w)$, and $\mathcal L(\tau(w))\subseteq\mathcal L(v)\cup\mathcal L(w)$. Therefore, $(\tau(v),\tau(w))$ is frugal and has smaller length than $(v,w)$.

Now let us assume that $(v,w)$ is not frugal, and hence disjoint. We see that then $X=\mathcal L(v)\sqcup\mathcal L(w)$. Assume without loss of generality that the multiplier $a$ of the Whitehead automorphism $\tau_1$ is in $\Lambda(v)$. Let the Whitehead automorphism $\tau$ be defined on the generator set $X$ as follows.
\begin{equation*}
\tau(x):=
\begin{cases}
\tau_1(x)&\text{if }x\in\mathcal L(v)\\
x&\text{if }x\in\mathcal L(w)\\
\end{cases}
\end{equation*}

By Lemma \ref{badmult}, since $a\notin\Lambda(w)$, we know that $\abs{\tau_1(w)}\ge\abs w$. Since $\abs{\tau_1(v)}+\abs{\tau_1(w)}<\abs v+\abs w$, we have that $\abs{\tau_1(v)}<\abs v$. It is clear that $\tau(v)=\tau_1(v)$ and $\tau(w)=w$. Therefore, $\abs{\tau(v)}+\abs{\tau(w)}<\abs v+\abs w$. Moreover, $\mathcal L(\tau(v))=\mathcal L(\tau_1(v))\subseteq\mathcal L(v)$ and $\mathcal L(\tau(w))=\mathcal L(w)$. Therefore, $(\tau(v),\tau(w))$ is disjoint.
\end{proof}
\end{lem}

\begin{cor}\label{down}
If an orbit under the action of $\Aut(F)$ contains a good pair, then it contains a good pair that has minimal length.
\begin{proof}
If $(v,w)$ is a good pair then we can repeatedly apply Whitehead autmorphisms to it such that the image is a good pair of smaller length, until we obtain a good pair that has minimal length.
\end{proof}
\end{cor}

We can now prove Theorem \ref{algcyc}.

\begin{thm}
It is decidable whether or not two conjugacy classes $v$ and $w$ or $F$ are both elliptic to some proper free product decomposition of $F$.

\begin{proof}
By Theorem \ref{orbit}, for any pair of cyclic words, we can determine whether or not they are in the same orbit as $(v,w)$ under the action of $\Aut(F)$. We can also determine the minimal length of a pair in the orbit of $(v,w)$. There are finitely many pairs of that length, so for each pair of that length we can check whether or not it is in the orbit of $(v,w)$ and whether or not it is good. Therefore, we can determine whether or not there exists a good pair in the orbit of $(v,w)$ that has minimal length.

If such a pair exists, then by Lemma \ref{auto}, there exists a proper free product decomposition $A*B$ of $F$ such that both $v$ and $w$ are elliptic to $A*B$. Conversely, if no such pair exists, then by the contrapositive of Corollary \ref{down}, the orbit of $(v,w)$ does not contain a good pair, so there does not exist a proper free product decomposition such that both $v$ and $w$ are elliptic to it.
\end{proof}
\end{thm}

This algorithm requires finding all pairs of minimal length in the orbit of $(v,w)$ and determining whether or not any of them are good. We will show that it suffices to check just one pair of minimal length.

\begin{lem}\label{across}
If $(v,w)$ is a good pair of minimal length and $\tau$ is a Whitehead automorphism such that $\bigl(\tau(v),\tau(w)\bigr)$ also has minimal length, then $\bigl(\tau(v),\tau(w)\bigr)$ is also good.

\begin{proof}
If $\tau$ is a relabeling automorphism then we have that $\Lambda(\tau(v))=\tau(\Lambda(v))$ and $\Lambda(\tau(w))=\tau(\Lambda(w))$, so $(\tau(v),\tau(w))$ is good. We will henceforth consider the case where $\tau$ is not a relabeling automorphism. Let the multiplier of $\tau$ be $a$.

Let $(v,w)$ be frugal. If $a\in\Lambda(v)\cup\Lambda(w)$, then we have $\mathcal L(\tau(v))\subseteq\mathcal L(v)\cup\mathcal L(w)$ and $\mathcal L(\tau(w))\subseteq\mathcal L(v)\cup\mathcal L(w)$, so $(\tau(v),\tau(w))$ is frugal. If, on the other hand, $a\notin\Lambda(v)\cup\Lambda(w)$, then by Lemma \ref{badmult}, $\abs{\tau(v)}\ge\abs v$ and $\abs{\tau(w)}\ge\abs w$. Since $\abs{\tau(v)}+\abs{\tau(w)}=\abs v+\abs w$, we know that $\abs{\tau(v)}=\abs v$ and $\abs{\tau(w)}=\abs w$. Again using Lemma \ref{badmult}, we conclude that $\tau(v)=v$ and $\tau(w)=w$, so $(\tau(v),\tau(w))$ is good.

If $(v,w)$ is not frugal, then it is disjoint, and so $X=\mathcal L(v)\sqcup\mathcal L(w)$. Without loss of generality, assume that $a\in\Lambda(v)$. Then $\mathcal L(\tau(v))\subseteq\mathcal L(v)$. Assume for contradiction that $\abs{\tau(v)}<\abs v$. Then define a Whitehead automorphism $\tau'$ on the generating set $X$ as follows.
\begin{equation*}
\tau'(x):=
\begin{cases}
\tau(x)&\text{if }x\in\mathcal L(v),\\
x&\text{if }x\in\mathcal L(w).\\
\end{cases}
\end{equation*}
We have that $\abs{\tau'(v)}<\abs v$ and $\abs{\tau'(w)}=\abs w$, so $\abs{\tau'(v)}+\abs{\tau'(w)}<\abs v+\abs w$, which contradicts the assumption that $(v,w)$ is of minimal length. Therefore, $\abs{\tau(v)}\ge\abs v$. By Lemma \ref{badmult}, $\abs{\tau(w)}\ge\abs w$. Since $\abs{\tau(v)}+\abs{\tau(w)}=\abs v+\abs w$, we conclude that $\abs{\tau(w)}=\abs w$. Again using Lemma \ref{badmult}, $\tau(w)=w$, so $\mathcal L(\tau(w))=\mathcal L(w)$, and so $(\tau(v),\tau(w))$ is disjoint.
\end{proof}
\end{lem}

\begin{prop}\label{mingood}
If an orbit under the action of $\Aut(F)$ contains a good pair, then every pair in the orbit of minimal length is good.

\begin{proof}
By Corollary \ref{down}, there exists a good pair $(v,w)$ of minimal length in the orbit. By Theorem \ref{white}, if $(v',w')$ is a pair in the orbit of $(v,w)$ of minimal length then we can get from $(v,w)$ to $(v',w')$ by composing a sequence of Whitehead automorphisms such that, after each automorphism, the image of $(v,w)$ remains of minimal length. Therefore, by Lemma \ref{across}, the image of $(v,w)$ is good after each automorphism, so $(v',w')$ is good.
\end{proof}
\end{prop}

Thus, we can check whether or not $v$ and $w$ are elliptic to a common proper free product decomposition by applying length decreasing Whitehead automorphisms until we arrive at a pair $(v',w')$ of minimal length. The words $v$ and $w$ are elliptic to a common proper free product decomposition if and only if the pair $(v',w')$ is good.

\section*{Acknowledgements}
I would like to thank Kim Whittlesey and Ilya Kapovich at the University of Illinois at Urbana-Champaign for introducing me to this topic. I would also like to thank Matthew Day at Caltech for his advice on revising this paper. Finally, I would like to thank the anonymous reviewers for their comments, particularly for bringing Stalling's paper \cite{swgoh} to my attention.

\bibliographystyle{amsplain}
\bibliography{Cite}

\end{document}